\documentclass[journal]{IEEEtran}

\usepackage{xcolor}
\usepackage{graphicx}
\usepackage{ifpdf}
\usepackage{epstopdf}
\usepackage{ifpdf}
\ifpdf
\DeclareGraphicsExtensions{.eps}
\else
\DeclareGraphicsExtensions{.eps}
\fi
\ifCLASSOPTIONcompsoc
\usepackage[caption=false,font=normalsize,labelfon
t=sf,textfont=sf]{subfig}
\else
\usepackage[caption=false,font=footnotesize]{subfi
	g}
\fi

\usepackage[cmex10]{amsmath}
\usepackage{amssymb}  
\usepackage{amsthm}
\usepackage[noadjust]{cite}
\usepackage{enumitem}
\usepackage{multirow}
\newcommand\MYhyperrefoptions{bookmarks=true,bookmarksnumbered=true,
	pdfpagemode={UseOutlines},plainpages=false,pdfpagelabels=true,
	colorlinks=true,linkcolor={black},citecolor={black},urlcolor={black},
	pdftitle={},
	pdfsubject={},
	pdfauthor={},
	pdfkeywords={}}

\ifCLASSINFOpdf
\usepackage[\MYhyperrefoptions,pdftex]{hyperref}
\usepackage{breakurl}
\else
\usepackage[\MYhyperrefoptions,breaklinks=true,dvips]{hyperref}
\usepackage{breakurl}
\fi

\usepackage{nomencl}
\makenomenclature

\usepackage{etoolbox}
\renewcommand\nomgroup[1]{%
	\item[\bfseries
	\ifstrequal{#1}{C}{Constants and parameters}{%
		\ifstrequal{#1}{I}{Sets and indices}{%
			\ifstrequal{#1}{F}{Functions}{%
				\ifstrequal{#1}{S}{Computed quantities}{%
			\ifstrequal{#1}{V}{Decision variables}{}}}}}%
	]}

\newcommand{\eq}[2]{
	\begin{IEEEeqnarray}{#1}
		#2
\end{IEEEeqnarray}}

\newcommand{\mr}[1]{\mathrm{#1}}

\newcommand{\mbb}[1]{\mathbb{#1}}

\newcommand{\mc}[1]{\mathcal{#1}}

\DeclareMathOperator*{\minimize}{minimize}
\DeclareMathOperator*{\subjectto}{subject~to}

\begin{document}

\title{Computationally Efficient Solutions for Large-Scale Security-Constrained Optimal Power Flow}

\author{Mohammadhafez~Bazrafshan,~Kyri~Baker,~\IEEEmembership{Member,~IEEE},~and Javad~Mohammadi,~\IEEEmembership{Senior Member,~IEEE}%
\thanks{\noindent   This work acknowledges support by the  U.S.  Department  of  Energy’s  Advanced Research Projects Agency-Energy (ARPA-E) under award DE-AR0001082.  This work utilized the Summit supercomputer, which is supported by the National Science Foundation (awards ACI-1532235 and ACI-1532236), the University of Colorado Boulder, and Colorado State University. 
The Summit supercomputer is a joint effort of the University of Colorado Boulder and Colorado State University.
}
}

\maketitle

\begin{abstract}
In this paper, we discuss our approach and algorithmic framework for solving large-scale security constrained optimal power flow (SCOPF) problems. SCOPF is a mixed integer non-convex optimization problem that aims to obtain the minimum dispatch cost while maintaining the system N-1 secure. Finding a feasible solution for this problem over large networks is challenging and this paper presents contingency selection, approximation methods, and decomposition techniques to address this challenge in a short period of time. The performance of the proposed methods are verified through large-scale synthetic and actual power networks in the Grid Optimization (GO) competition organized by the U.S. Advanced Research Projects Agency-Energy (ARPA-E). As many prior works focus on small-scale systems and are not benchmarked using validated, publicly available datasets, we aim to present a practical solution to SCOPF that has been proven to achieve good performance on realistically sized (30,000 buses) networks. 


\end{abstract}

\begin{IEEEkeywords}
Security-constrained optimal power flow, Optimization decomposition, High performance computing \vspace{-3mm}
\end{IEEEkeywords}

\printnomenclature

\vspace{-2.5mm}
\section{Introduction}
The security constrained optimal power flow (SCOPF) problem seeks to achieve economically optimal steady-state setpoints for power system controllers to operate under normal operation and during contingencies, e.g., loss of a generator, transformer, or transmission line. This extends the traditional optimal power flow (OPF) problem to adhere to the N-1 security criteria. Put differently, this problem ensures equity of load and generation 
while conforming to network constraints (e.g. voltage and generation limits) in the event of a single element failure (e.g., line, generator, or transformer).

This paper discusses our optimization framework which solves the
SCOPF problems using the formulation and grid configurations developed by the Department of Energy's the Advanced Research Projects Agency-
Energy (ARPA-E) as part of the Grid Optimization (GO) competition. Our proposed solution was ranked as one of the top performing solutions in the competition on large-scale, realistic networks (up to 30,000 buses), and was ranked second-place on actual industry power networks. In this paper we discuss how these results were achieved and present the performance of different solutions that we examined and tested as part of our participation in the GO competition.

The SCOPF problem formulations follows that of a two-stage non-convex stochastic optimization problem where the stochasticity is, in essence, deterministically represented by equally probable scenarios,  see e.g.,~\cite{Alsac1974}. 
In this setup, the first stage formulation, also known as the preventive stage, is associated with normal operation of the power grid. The second stage, also known as the corrective stage represents operation under contingencies. The two main bottle-necks in solving the SCOPF problem are (i) the non-convexity of the power flow equations and (ii) the computational complexity arising from contingency inclusion posed as second-stage constraints coupled with first-stage normal operation constraints.

Addressing these two challenges have been vastly studied in the literature. 
To cope with the non-convexity challenge, researchers have studied a wide-rage of convex modifications of the power flow equations, e.g., the linear program (LP) in~\cite{mohammadi2017fully}, the second order cone program in~\cite{Wu2018}, and the semidefinite program in~\cite{Madani2016}. To address the computational complexity issue, researchers have proposed contingency filtering and screening to reduce the set of considered contingencies. These techniques aim to find the most critical failures based on their impact on the system dispatch, e.g.~\cite{Ardakani2013,Madani2017}. To further relieve the computational burden arising from contingency constraints, decomposition techniques such as Benders~\cite{li2008decomposed,mohammadi2013benders}, Lagrangian methods~\cite{biskas2004decentralised,kim2001method}, and consensus-based methods \cite{Mohammadi2018} are suggested in the literature. Prior methods struggle to solve the comprehensive SCOPF formulation (implemented in the industry software that is used by power grid operators) and often make simplifying assumptions to make the original problem tractable. Translating theoretical methodologies to industry requirements is, however, not straightforward. In this paper we use the ARPA-E GO competition's formulation and data-sets which are adopted from realistic formulation used by major Independent System Operators (ISOs).

For instance, in a real-time setting, the preventive stage of the SCOPF must produce a feasible base-case solution for networks with as many as 30,000 buses under 10 minutes~\cite{ARPAEScore}. In such a time-restricted setting, feasibility must be prioritized over optimality. Furthermore, to examine the quality of the preventive solution, tens of thousands of contingencies need to be considered and solved in the second stage of the SCOPF problem. Hence a comparison between optimality and feasibility is also necessary for corrective decisions.

In this paper we propose computationally efficient solutions for ultra large-scale power grids.
To this end, we will present two contingency selection procedures, three formulation approximations combined with algorithmic decomposition that we used for handling many equations and variables associated with large-scale SCOPF problem. Implementing various combinations of these techniques yield twelve distinct approaches. Impartial comparisons between these algorithms based on the trade-off between their optimality and time efficiency guided our top-performing solution approach. We highlight here that although our methods are discussed in the context of solving ARPA-E GO's formulation~\cite{ARPAEFormulation}, they can be adopted to solve similar representations of the SCOPF problem. This formulation follows the standard industry practice and is comprised of preventive or pre-contingency (hereafter named code 1) and post-contingency corrective actions (hereafter named code 2). 

This paper is organized as follows. In Section~\ref{sec:scopf}, the mathematical formulation of the SCOPF problem--inspired by the GO competition---is presented. Submission of the solutions and their evaluation are also briefly outlined. Our proposed solutions are detailed in Sections~\ref{sec:code1heuristics} and~\ref{sec:code2heuristics} followed by discussions on their implementation in Section~\ref{sec:implementation}. Computational results are detailed in  Section~\ref{sec:numericaltests}. Finally, section~\ref{sec:conclusion} concludes this paper.

\vspace{-3mm}
\section{SCOPF formulation and evaluation}
\label{sec:scopf}
\subsection{Problem formulation} 
\label{sec:scopfproblem}
We give a high level description of the SCOPF problem as presented by the ARPA-E GO competition formulation~\cite{ARPAEFormulation}. The problem comprises a \emph{base case} denoted by the index  $k=0$ and a number of $K$ \emph{contingency cases} denoted by indices  $k$ in the contingency set $\mc{K}=\{1,\ldots,K\}$. A contingency refers to an outage of one component from the set of generators, transmission lines, and transformers that are online in the  base case.  
\nomenclature[C]{$K$}{Number of contingencies}
\nomenclature[I]{$k$}{Contingency index, $k=0$ denotes  base case}
\nomenclature[I]{$\mc{K}$}{Contingency set}
\emph{Primary decision variables} include voltage magnitudes, voltage angles, real and reactive power injections from generators and controllable shunt susceptances.  For contingencies $k \in \mc{K}$, an additional post-contingency  real power adjustment variable $\Delta_k$  is also included in the primary decision vector.  Let us denote the number of primary decision variables in the base and contingency cases by $n_{x_k}$   and collect the primary variables   in vectors $x_k \in \mbb{R}^{n_{x_k}}$ for  $k=\{0\} \cup \mc{K}$.  

\nomenclature[V]{$\Delta_k$}{Generator real power adjustment for contingency $k$}
\nomenclature[C]{$n_{x_k}$}{Number of primary decision variables in the base case} 
\nomenclature[V]{$x_k$}{Primary decision variables}
\emph{Slack decision variables} 
are introduced in the comprehensive SCOPF formulation~\cite{ARPAEFormulation} to allow for soft constraints; i.e., those constraints that may have violations. These slack variables are penalized in the objective of the problem.
Slack variables are denoted by $s_k \in \mbb{R}^{n_{s_k}}$  where $n_{s_k}$ denotes the number of slack variables for $k=\{0\} \cup \mc{K}$.  The vector  $s_k$ is further compartmentalized to $s_k=[s_k^{f+'}, s_k^{f-'}, s_k^{g'}]'$ where $s_k^{f+}$, $s_k^{f-}$ are slack variables for soft power flow constraints, $s_k^{g}$ are slack variables for soft transmission line and transformer rating constraints, and $'$ denotes transposition. 
\nomenclature[V]{$s_k$}{Slack variables per index $k$}
\nomenclature[C]{$n_{s_k}$}{Number of slack variables per index $k$}
\nomenclature[V]{$s_k^{f+},$}{Slack for positive violation of  power flow constraints per index $k$}
\nomenclature[V]{$s_k^{f-},$}{Slack for negative violation of power flow constraints per index $k$}
\nomenclature[V]{$s_k^{g}$}{Slack for violation of line/transformer rating per index $k$}

\noindent We thus consider the following SCOPF formulation: 
\eq{rCl}{\label{eqgroup:scopf}\IEEEyesnumber \IEEEyessubnumber* 
	\minimize_{x_0,s_0,x_k,s_k} \: c(x_0)+c_0(s_0)&+&\frac{1}{K} \sum\limits_{k\in \mc{K}} c_k(s_k) \IEEEeqnarraynumspace \label{eq:scopfobj} \\
	\subjectto \: f_0(x_0)&=&s_{0}^{f+}-s_{0}^{f-} \hfill \label{eq:scopff0} \\
	g_0(x_0)&\le& s_0^{g}\label{eq:scopfg0}\\
	f_k(x_k) &=& s_{k}^{f+}-s_{k}^{f-}, \: k \in \mc{K} \label{eq:scopffk}\\
	g_k(x_k)&\le&s_{k}^{g} \: k \in \mc{K}\label{eq:scopfgk}\\
	\IEEEeqnarraymulticol{3}{c}{h_k(x_0^{\mr{p},\mr{q},\mr{v}},x_k^{\mr{p},\mr{q},\mr{v},\Delta}) =0, \: k \in \mc{K}} \label{eq:scopf0k}\\
	\IEEEeqnarraymulticol{3}{c}{x_k \in \mc{X}_k, s_k \in \mc{S}_k, \: k \in \{0\} \cup \mc{K}.\IEEEeqnarraynumspace \label{eq:scopfboxes}}
}

 The objective in~\eqref{eq:scopfobj} comprises piece-wise linear functions $c(\cdot)$, $c_0(\cdot)$, and $c_{k}(\cdot)$ which respectively calculate generator cost, penalty costs for the base case, and penalty costs for the contingency cases \cite[eqs (1)--(31)]{ARPAEFormulation}. Constraints~\eqref{eq:scopff0} and~\eqref{eq:scopfg0}  respectively account for soft power balance and  transmission line/transformer rating constraints for the base case.   Constraints~\eqref{eq:scopffk} and~\eqref{eq:scopfgk} respectively account for soft power balance and  transmission line/transformer rating constraints during contingencies.  Constraint~\eqref{eq:scopf0k} is the generator responses to contingencies. The superscripts $\mr{p}$,$\mr{q}$,$\mr{v}$, and $\Delta$ refer respectively to those indices of the  vector that correspond to generators' real power, reactive power,  voltages and real power adjustment. Lastly, the sets $\mc{X}_k$ and $\mc{S}_k$ are subsets of the extended real coordinate space in, respectively, $\mbb{R}^{n_{x_k}}$ and $\mbb{R}^{n_{s_k}}$. Hence, constraints expressed in~\eqref{eq:scopfboxes} express lower and upper bounds on primary and slack variables.
\nomenclature[F]{$f_k(\cdot)$}{Power balance equations with implicit power flows per index $k$}
\nomenclature[F]{$g_k(\cdot)$}{Line/transformer rating constraint with implicit power flows  per index $k$}
\nomenclature[F]{$h_k(\cdot)$}{Generator PV/PQ switching response equations per index $k$}
\nomenclature[F]{$c(\cdot)$}{Generator cost for base case}
\nomenclature[F]{$c_k(\cdot)$}{Penalty cost per index $k$}
\nomenclature[V]{$x^{\mr{p}}$}{Portion of $x$ pertaining to generators' real power}
\nomenclature[V]{$x^{\mr{q}}$}{Portion of $x$ pertaining to generators' reactive power}
\nomenclature[V]{$x^{\mr{v}}$}{Portion of $x$ pertaining to generators' voltages}
\nomenclature[V]{$x^{\mr{\Delta}}$}{Portion of $x$ pertaining to generators' real power adjustment}
\nomenclature[I]{$\mc{X}_k$}{Rectangular feasible region in $\mbb{R}^{n_{x_k}}$ for primary decision variables $x_k$}
\nomenclature[I]{$\mc{S}_k$}{Rectangular feasible region in $\mbb{R}^{n_{s_k}}$ for slack decision variables $s_k$}

Although we do not delve into specifics of every function in the objective or constraints throughout this manuscript, it is useful to highlight the  challenges in solving~\eqref{eqgroup:scopf}.  The first challenge is the nonlinearity and nonconvexity of  soft power balance constraints~\eqref{eq:scopff0} and~\eqref{eq:scopffk}. The second  is handling generator response constraints~\eqref{eq:scopf0k} which are nonlinear, nonconvex,  and nondifferentiable. These constraints may be represented by min/max operators or by a disjunction of linear constraints. Inputting ~\eqref{eq:scopf0k} into optimization solvers require introduction of binary variables~\cite[Sections 3.14 and 3.15]{ARPAEFormulation}. 

Last but not least, a third challenge is the size of the problem. Indeed, this problem is a two-stage stochastic program where the base and contingency variables comprise the first and second stage decisions, respectively. As the number of buses and branches in the networks increases, the number of contingencies $K$ also inevitably rises. Networks with thousands of buses, branches, and thousands of contingencies are common in the GO competition  datasets~\cite{ARPAEDatasets}.  The coupling introduced in constraint~\eqref{eq:scopf0k} prevents a straightforward parallelization of the solution methodology rendering the problem intractable for larger networks.  


\subsection{Solution submission procedure}
 Computed base case primary variables, $x_0^*$, are printed in \texttt{solution1.txt} and computed contingency primary variables,  $x_k^*$ for $k \in \mc{K}$, are printed in \texttt{solution2.txt}.   Slack variables are not submitted for they are recomputed  using the evaluation procedure~\cite[Appendix F]{ARPAEFormulation}.
Accordingly, the algorithm  has two main parts:  Code 1 and code 2, respectively   outputting \texttt{solution1.txt} and \texttt{solution2.txt}. 
\nomenclature[S]{$x_k^*$}{Computed values for primary decision vector $x_k$}

\vspace{-2mm}


\vspace{-2mm}
\subsection{Evaluation procedure}\label{sec:evaluation}

The performance evaluation criteria used to measure the performance our proposed solutions in described here according to   ARPA-E GO's evaluation metrics \cite{ARPAEScore}.
Prior to elaborating on the evaluation procedure, an explanation is in order to construct feasible solutions $\breve{x}_0$ and $\breve{x}_k$ for $k \in \mc{K}$ that yield the worst-case cost for a power grid:  
\nomenclature[S]{$\breve{x}_k$}{Values for primary decision vector $x_k$ that yield the worst-case cost}
\begin{enumerate}
	\item Set the variables in indices of $\breve{x}_0$ and $\breve{x}_k$ corresponding to the  voltage magnitudes, real and reactive power generation to the midpoint of their bound constraints described by the sets $\mc{X}_k$ for $k \in \{0\} \cup \mc{K}$.
	\item  For every generator contingency $k$, re-set those variables of $\breve{x}_k$ corresponding to the real and reactive power outputs at contingency $k$ to zero.
	\item Set the variables in indices of $\breve{x}_0$ and $\breve{x}_k$ corresponding to the voltage angles, controllable shunt susceptances, and real power generation adjustments to zero. 
\end{enumerate} 

Based on the form of the constraint~\eqref{eq:scopf0k} (for instance, \cite[eqs. (85) and (93)]{ARPAEFormulation}), it is straightforward to see that the worst-case solution procedure ensures that $h_k(\breve{x}_0^{\mr{p},\mr{q},\mr{v}},\breve{x}_k^{\mr{p},\mr{q},\mr{v},\Delta}) =0$.
 Plugging in $\breve{x}_k$ in $f_k(\cdot)$ and $g_k(\cdot)$ yields $\breve{s}_k$.  The worst-case cost (or sometimes referred to as  score) is then given by  
 \eq{rCl}{c^{\mr{slack}}=c(\breve{x}_0)+c_0(\breve{s}_0)+\frac{1}{K} \sum\limits_{k\in \mc{K}} c_k(\breve{s}_k). \label{eq:cslack}}
\nomenclature[S]{$c^{\mr{slack}}$}{worst-case cost (worst score)}
We now describe the following evaluation procedure: 

\begin{enumerate}
	\item After the specified amount of time for code 1, read $x_0^*$ from file \texttt{solution1.txt}. If the file does not exist or any other errors arise (such as  incorrect formatting) declare the worst case score and terminate.
	\item After the specified amount of time for code 2, read $x_k^*$ for $k \in \mc{K}$  from file \texttt{solution2.txt}. If the file does not exist or any other errors arise declare the worst case score and terminate.
	\item After successfully reading $x_0^*$ and $x_k^*$ for $k \in \mc{K}$, if $h_k(x_0^{{\mr{p},\mr{q},\mr{v}} *},x_k^{{\mr{p},\mr{q},\mr{v},\Delta} *}) \neq 0$ for any $k \in \mc{K}$, or if $x_k \not\in \mc{X}_k$ for any $k \in \{0\} \cup \mc{K}$, declare the worst case score and terminate.
	\item If $h_k(x_0^{\mr{p},\mr{q},\mr{v}*},x_k^{\mr{p},\mr{q},\mr{v},\Delta *}) =0$ for all $k \in\mc{K}$ and if $x_k^* \in \mc{X}_k$ for all $k \in \{0\} \cup \mc{K}$, plug in $x_k^*$ in $f_k(.)$ and $g_k(.)$ to yield $s_k^*$ for $k \in \{0\} \cup \mc{K}$. Set the feasible score as
	 \eq{rCl}{c^{*}=c(x_0^*)+c_0(s_0^*)+\frac{1}{K} \sum\limits_{k\in \mc{K}} c_k(s_k^*). \label{eq:cfeasible}}
	 \item Declare network score as $c^{\mr{score}}=\min\{c^{*}, c^{\mr{slack}}\}$. 
	\end{enumerate}
\nomenclature[S]{$c^*$}{Feasible score}
\nomenclature[S]{$c^{\mr{score}}$}{Network score}
Having reviewed the problem formulation, implementation, and evaluation procedures, we now discuss our code 1 (representing preventative SCOPF) and code 2 (representing corrective SCOPF) solutions. 

\vspace{-3mm}
\section{Preventative SCOPF (``code 1") solution}
\label{sec:code1heuristics}
\vspace{-1mm}
In this section, we present our algorithm for finding a solution for code 1 in \emph{real time}, the restriction that code 1 must produce \texttt{solution1.txt} within  ten minutes of wall-clock time.  Our \emph{offline} algorithm is based on similar principles with slight improvements.  In this section, we discuss several methods that we used for developing code 1.  These solutions are described next in the three categories of contingency selection, formulation approximation, and decomposition approaches.

\vspace{-3mm}
\subsection{Contingency selection}
\label{sec:conselection}
 This subsection presents two practical contingency selection methods.  The goal  is to trim down the set $\mc{K}$ with cardinality $K$ to a   set $\hat{\mc{K}}$ with cardinality $\hat{K} < K$ so that the SCOPF problem~\eqref{eqgroup:scopf} becomes tractable.  The first method is based on nominal sizes or \emph{ratings} of generators, transmission lines, and transformers. The second is based on knowledge of a base case operating point from the data files--referred to as contingency selection based on \emph{real-time} data.
 \nomenclature[I]{$\hat{\mc{K}}$}{Reduced contingency set}
 \nomenclature[C]{$\hat{K}$}{Number of contingencies in reduced contingency set}
 
 \subsubsection{Contingency selection based on asset ratings}
Contingency selection based on ratings is  as follows:

 \begin{enumerate}[label=(\roman*)]
 	\item Select an integer $\hat{K} < K$.
 	\item For every contingency $k \in \mc{K}$ that describes the outage of a generator, let  $\bar{p}_k$  be this out-of-service generator's \emph{base-case} real  power upper bound obtained from $\mc{X}_0$. The \emph{rank} of such a  contingency is calculated as
 	\eq{rCl}{\mr{rank}_{k}=\bar{p}_k}
 	\item For every contingency $k \in \mc{K}$ that describes the outage of a transmission line or transformer, let  $\bar{g}_k$ be this out-of-service branch's  \emph{base-case} maximum apparent power rating.
 	The \emph{rank} of such a  contingency is calculated as
 	\eq{rCl}{\mr{rank}_{k}=\bar{g}_k.}
 	\item Collect the computed ranks $\mr{rank}_k$ in a vector and use a  descending sort to obtain the vector $\mr{rank}$.   Output the set $\hat{\mc{K}} \subset \mc{K}$ that collects  indices $\hat{k} \in \mc{K}$ where  $\mr{rank}_{\hat{k}}$ is among the first $\hat{K}$ entries of $\mr{rank}$. 
 \end{enumerate}
\nomenclature[C]{$\bar{p}_k$}{Base-case real power upper bound of generator out of service in contingency $k$}
\nomenclature[C]{$\bar{g}_k$}{Base-case apparent  rating of line/transformer out of service in contingency $k$}
\nomenclature[S]{$\mr{rank}_k$}{Rank of a contingency $k$ indicating its importance in optimization}
\nomenclature[I]{$\hat{k}$}{Index for reduced contingency set}
\subsubsection{Contingency selection based on real-time computations}
Suppose a previously calculated base case operating point solution is given by $\tilde{x}_0$. The algorithm for contingency selection based on real-time computations is as follows:
 \begin{enumerate}[label=(\roman*)]
	\item Select an integer $\hat{K} < K$.
	\item For every contingency $k \in \mc{K}$ that describes the outage of a generator, let  $\tilde{p}_k$ and $\tilde{q}_k$  be this out-of-service generator's real  power and reactive power generation from the given operating point $\tilde{x}_0$. The \emph{rank} of such a  contingency is calculated as
	\eq{rCl}{\mr{rank}_{k}=\sqrt{\tilde{p}_k^2+\tilde{q}_k^2}.}
	\item For every contingency $k \in \mc{K}$ that describes the outage of a branch (transmission line or transformer), let $\tilde{p}_{k,o}$ and $\tilde{q}_{k,o}$ together with $\tilde{p}_{k,d}$ and $\tilde{q}_{k,d}$ collect respectively the real and reactive power flows at the origin and destination of the corresponding branch computed from $\tilde{x}_0$. The \emph{rank} of such a  contingency is calculated as
	\eq{rCl}{\mr{rank}_{k}=0.5\sqrt{\tilde{p}_{k,o}^2+\tilde{q}_{k,o}^2}+0.5\sqrt{\tilde{p}_{k,d}^2+\tilde{q}_{k,d}^2}.}
	\item Collect the computed ranks $\mr{rank}_k$ in a vector and use a  descending sort to obtain the vector $\mr{rank}$.   Output the set $\hat{\mc{K}} \subset \mc{K}$ that collects  indices $\hat{k} \in \mc{K}$ where  $\mr{rank}_{\hat{k}}$ is among the first $\hat{K}$ entries of $\mr{rank}$. 
\end{enumerate}
\nomenclature[C]{$\tilde{p}_k$}{Last-known  operating point of real power for generator out of service in contingency $k$}
\nomenclature[C]{$\tilde{q}_k$}{Last-known operating point of reactive power for  of generator out of service in contingency $k$}
\nomenclature[C]{$\tilde{x}_k$}{Last-known operating point of the network}
\nomenclature[C]{$\tilde{p}_{k,o}$}{Last-known operating point of  real power flow from origin of branch out of service in contingency $k$}
\nomenclature[C]{$\tilde{q}_{k,o}$}{Last-known operating point of  reactive power flow from origin of branch out of service in contingency $k$}
\nomenclature[C]{$\tilde{p}_{k,d}$}{Last-known operating point of  real power flow from destination of branch out of service in contingency $k$}
\nomenclature[C]{$\tilde{q}_{k,d}$}{Last-known operating point of  reactive power flow from destination of branch out of service in contingency $k$}

\vspace{-2mm}
\subsection{Approximate Formulations}
\label{sec:approxformulations}
Building on our contingency selection methods, this section propose approximate solutions for solving the \emph{reduced} SCOPF problem~\eqref{eqgroup:scopf} obtained by replacing the \emph{reduced} set $\hat{\mc{K}}$   for $\mc{K}$.
Although the proposed contingency selection methods reduce the size of contingency set, solving the reduced SCOPF problem for large-scale network still faces two main difficulties. These two challenges originate from 
the nonconvexity of constraints~\eqref{eq:scopff0} and~\eqref{eq:scopffk} as well as the nonconvexity and nondifferentiability of  function $h(\cdot)$ in constraint~\eqref{eq:scopf0k}.  Depending on various handling of the aformentioned constraints, three  approximations of the \emph{reduced} problem are given next.

\subsubsection{NBNC formulation}
The NBNC formulation is given as
\eq{rCl}{\label{eqgroup:nbnc}\IEEEyesnumber \IEEEyessubnumber* 
	\minimize_{x_0,s_0,x_k,s_k} \: c(x_0)&+&c_0(s_0)+\frac{1}{\hat{K}} \sum\limits_{k\in \hat{\mc{K}}} c_k(s_k) \IEEEeqnarraynumspace \label{eq:nbncobj} \\
	\subjectto \: f_0(x_0)&=&s_{0}^{f+}-s_{0}^{f-} \hfill \label{eq:nbncf0} \\
	g_0(x_0)&\le& s_0^{g}\label{eq:nbncg0}\\
	f_k(x_k) &=& s_{k}^{f+}-s_{k}^{f-}, \: k \in \hat{\mc{K}} \label{eq:nbncfk}\\
	g_k(x_k)&\le&s_{k}^{g}, \: k \in \hat{\mc{K}}\label{eq:nbncgk}\\
	\IEEEeqnarraymulticol{3}{c}{\hbar_k(x_0^{\mr{p},\mr{v}},x_k^{\mr{p},\mr{v},\Delta}) =0, \: k \in \hat{\mc{K}} \label{eq:nbnc0k}}\\
	\IEEEeqnarraymulticol{3}{c}{x_k \in \mc{X}_k, s_k \in \mc{S}_k, \: k \in \{0\} \cup \hat{\mc{K}}. \label{eq:nbncboxes}}
}
 Here, power balance and branch rating constraints retain their nonlinear, nonconvex form (NBNC stands for  Nonlinear Base Nonlinear Contingencies).  However, constraint~\eqref{eq:nbnc0k} is introduced  where a linear function $\hbar_k(\cdot)$ is used in place of $h_k(\cdot)$ to provide a linear inner approximation of constraint ~\eqref{eq:scopf0k}. The function $\hbar_k(\cdot)$ ignores the dependence of~\eqref{eq:scopf0k} on reactive powers.
Since all the base case variables  $x_0$ of the original SCOPF problem~\eqref{eqgroup:scopf} are present in this formulation, the solution $x_0^*$  is feasible for the base case and can be output as \texttt{solution1.txt}.  Recall that code 1 is responsible only for providing a feasible solution for the base case.
\nomenclature[F]{$\hbar(\cdot)$}{Linear inner approximation of $h_k(\cdot)$}
\subsubsection{NBLC formulation}
This formulation is presented as
\eq{rCl}{\label{eqgroup:nblc}\IEEEyesnumber \IEEEyessubnumber* 
	\minimize_{x_0,s_0,\hat{x}_k,\hat{s}_k} \: c(x_0)&+&c_0(s_0)+\frac{1}{\hat{K}} \sum\limits_{k\in \hat{\mc{K}}} \hat{c}_k(\hat{s}_k) \IEEEeqnarraynumspace \label{eq:nblcobj} \\
	\subjectto \: f_0(x_0)&=&s_{0}^{f+}-s_{0}^{f-} \hfill \label{eq:nblcf0} \\
	g_0(x_0)&\le& s_0^{g}\label{eq:nblcg0}\\
	\hat{f}_k(\hat{x}_k) &=& \hat{s}_{k}^{f+}-\hat{s}_{k}^{f-}, \: k \in \hat{\mc{K}} \label{eq:nblcfk}\\
		\IEEEeqnarraymulticol{3}{c}{\hat{h}_k(x_0^{\mr{p}},\hat{x}_k^{\mr{p},\Delta}) =0, \: k \in \hat{\mc{K}} \label{eq:nblc0k}}\\
	\IEEEeqnarraymulticol{3}{c}{x_0 \in\mc{X}_0, s_0\in \mc{S}_0, \hat{x}_k \in \hat{\mc{X}}_k, \hat{s}_k \in \hat{\mc{S}}_k, \: k \in \hat{\mc{K}}. \label{eq:nblcboxes}}
}
In NBLC (i.e., Nonlinear Base Linear Contingency), the power balance and branch ratings of the base case retain their original nonconvex form. For contingencies, power balance constraints are linearized using similar principles as those given by the DC power flow formulation; see, e.g.~\cite{Taylor2015} and references therein. The branch  ratings per contingency $k$ are ignored. In the NBLC formulation, due to the omission of reactive powers and voltages in the contingencies, the linear function $\hbar_k(\cdot)$ of~\eqref{eq:nbnc0k} is relaxed to $\hat{h}_k(\cdot)$ of~\eqref{eq:nblc0k} to serve as a linear surrogate of~\eqref{eq:scopf0k}.  Symbols $\hat{c}_k$, $\hat{s}_k$, $\hat{f}_k$, $\hat{x}_k$, $\hat{\mc{X}}_k$, and $\hat{\mc{S}}_k$ are the trimmed versions of the corresponding symbols per requirements of the DC approximation.
\nomenclature[F]{$\hat{h}_k(\cdot)$}{Relaxed version of $\hbar_k{\cdot}$ with only real power generation and its adjustment}
\nomenclature[F]{$\hat{f}_k(\cdot)$}{DC power balance equations with implicit power flows}
\nomenclature[F]{$\hat{c}_k(\cdot)$}{Trimmed version of $c_k(\cdot)$ with variables of DC power balance only}
\nomenclature[F]{$\hat{s}_k$}{Trimmed version of $s_k$ with variables of DC power balance only}
\nomenclature[F]{$\hat{x}_k$}{Trimmed version of $x_k$ with variables of DC power balance only}
\nomenclature[F]{$\hat{S}_k$}{Trimmed version of $\mc{S}_k$ with variables of DC power balance only}
\nomenclature[F]{$\hat{X}_k$}{Trimmed version of $\mc{X}_k$ with variables of DC power balance only}
Similar to the NBNC formulation, since all the base case variables  $x_0$ of the original SCOPF problem~\eqref{eqgroup:scopf} are present in~\eqref{eqgroup:nblc}, the  solution ${x}_0^*$  is feasible for the base case and may be used for \texttt{solution1.txt}.

\subsubsection{LBLC formulation}
The LBLC formulation is given as
\eq{rCl}{\label{eqgroup:lblc}\IEEEyesnumber \IEEEyessubnumber* 
	\minimize_{\hat{x}_0,\hat{s}_0,\hat{x}_k,\hat{s}_k} \: \hat{c}(\hat{x}_0)&+&\hat{c}_0(\hat{s}_0)+\frac{1}{\hat{K}} \sum\limits_{k\in \hat{\mc{K}}} \hat{c}_k(\hat{s}_k) \IEEEeqnarraynumspace \label{eq:lblcobj} \\
	\subjectto \: f_0(\hat{x}_0)&=&\hat{s}_{0}^{f+}-\hat{s}_{0}^{f-} \hfill \label{eq:lblcf0} \\
	f_k(\hat{x}_k) &=& \hat{s}_{k}^{f+}-\hat{s}_{k}^{f-}, \: k \in \hat{\mc{K}} \label{eq:lblcfk}\\
	\IEEEeqnarraymulticol{3}{c}{\hat{h}_k(\hat{x}_0^{\mr{p}},\hat{x}_k^{\mr{p},\Delta}) =0, \: k \in \hat{\mc{K}} \label{eq:lblc0k}}\\
	\IEEEeqnarraymulticol{3}{c}{\hat{x}_k \in \mc{X}_k, \hat{s}_k \in \hat{\mc{S}}_k, \: k \in \{0\} \cup \hat{\mc{K}}. \label{eq:lblcboxes}}
}
In the LBLC (stands for Linear Base Linear Contingency) formulation, power balance constraints are approximated using similar principles similar to that of DC power flows. Voltages, reactive powers,  and rating constraints are ignored. 
\nomenclature[F]{$\hat{c}(\cdot)$}{Trimmed version of $c(\cdot)$ with variables of DC power balance only}
Unlike  the NBNC and NBLC formulations, the base case variable $\hat{x}_0$ is a subcomponent of the the decision variable $x_0$. To elaborate, in order to calculate \texttt{solution1.txt}, the reactive powers, shunt susceptances, and voltage magnitudes remain to be calculated even after having solved~\eqref{eqgroup:lblc}. Although setting these remaining variables to midpoints of their corresponding regions technically yields a feasible base case solution, it is  desired to  achieve a feasible base case point with minimal slack penalties. Therefore, the following \emph{recovery} problem is defined to be solved immediately after solving LBLC:
\eq{rCl}{\label{eqgroup:recovery}\IEEEyesnumber \IEEEyessubnumber* 
	\minimize_{x_0,s_0} \: \|x_0^{\mr{p}}&-&\hat{x}_0^{\mr{p}*}\|_{2}+c_0(s_0) \hfill \IEEEeqnarraynumspace \label{eq:recoveryobj} \\
	\subjectto \: f_0(x_0)&=&s_{0}^{f+}-s_{0}^{f-} \hfill  \IEEEeqnarraynumspace\label{eq:recoveryf0} \\
	g_0(x_0)&\le& s_0^{g}\label{eq:recoveryg0}\\
	\IEEEeqnarraymulticol{3}{c}{x_0 \in \mc{X}_0, s_0 \in \mc{S}_0 \label{eq:recoveryboxes}
}}
In~\eqref{eqgroup:recovery}, $\hat{x}_0^{\mr{p}*}$ is the generator real power portion of the optimal  solution to~\eqref{eqgroup:lblc}. The goal of the recovery problem~\eqref{eqgroup:recovery} is to find a base case solution $x_0^*$ that is \emph{close} to  $\hat{x}_0^*$, the optimal solution of the LBLC problem~\eqref{eqgroup:lblc}---as measured by the Euclidean distance of  their subcomponents that collect generator real power outputs.

The NBNC, NBLC, and LBLC formulations offer approximate and more tractable versions of the original problem.
The LBLC formulation can be solved by CPLEX~\cite{CPLEX}. The NBNC, NBLC, and recovery  formulations can be fed to off-the-bench nonlinear programming (NLP) solvers such as SNOPT~\cite{snopt} and IPOPT~\cite{ipopt}. Solving these  extensive-form formluations using nonlinear optimization solvers are useful for assessing the performance of our proposed approximations and indeed the solutions provided by these have proved useful for smaller-sized networks.  However, for networks with thousands of buses and contingencies, extensive form approaches may fail to meet our \textit{real-time} criteria  laid out at the beginning of this section. 
We explore Benders decomposition~\cite{Geoffrion1972} as another technique to ensure achiving the best \texttt{solution1.txt} within the time limit of the \textit{real-time} criteria.

\vspace{-3mm}
\subsection{Benders decomposition}\label{sec:Benders}
\vspace{-1mm}
Benders decomposition is a widely used tool for disassembling the SCOPF problem into a series of base and contingency subproblems; see, e.g., \cite{Monticelli1987} and~\cite{Phan2014} for relevant examples. In what follows, we present modification of the NBNC, NBLC, and LBLC formulations to make them amenable to Benders decomposition. Applying the decomposition is then straightforward using instructions laid out in \cite{Conejo2006}. 
\subsubsection{Reformulating the NBNC problem~\eqref{eqgroup:nbnc} for decomposition}
 Per contingency $k \in \hat{\mc{K}}$, introduce new decision variables $y_k$ with the same dimension as $x_0^{\mr{p},\mr{v}}$. Remove constraint~\eqref{eq:nbnc0k}.  Insert the following constraints:
	\eq{rCl}{\IEEEyesnumber\label{eqgroup:nbncbenders} \IEEEyessubnumber*
	x_0^{\mr{p},\mr{v}}&=&y_k, \:  k \in \hat{\mc{K}} \label{eq:nbncbenderseq}\\
	\IEEEeqnarraymulticol{3}{c}{\hbar_k(y_k,x_k^{\mr{p},\mr{v},\Delta}) =0, \: k \in \hat{\mc{K}},} \label{eq:nbncbenders0k}}
rendering a new problem which is amenable to nonlinear Benders decomposition. The advantage of this technique is two-folds: (a) the subproblems can now be solved using $\hat{K}$ parallel processes and and (b) at every iteration $\ell$ of the Benders algorithm the incumbent solution $x_0(\ell)$ is a feasible base case solution which can be periodically output as \texttt{solution1.txt}.  This incumbent solution can serve as a back-up solution if the algorithm fails to converge in the specified time-limit.

\nomenclature[V]{$y_k$}{New decision variable introduced per index $k$ for Bender's decomposition}
\nomenclature[I]{$\ell$}{Iteration of Bender's algorithm}
\subsubsection{Reformulating the NBLC problem~\eqref{eqgroup:nblc} for decomposition}
Per contingency $k \in \hat{\mc{K}}$, introduce new decision variables $y_k$ with the same dimension as $x_0^{\mr{p}}$. Remove constraint~\eqref{eq:nblc0k}.  Insert the following constraints:
\eq{rCl}{\IEEEyesnumber \label{eqgroup:nblcbenders} \IEEEyessubnumber*
	x_0^{\mr{p}}&=&y_k, \:  k \in \hat{\mc{K}} \label{eq:nblcbenderseq}\\
	\IEEEeqnarraymulticol{3}{c}{\hat{h}_k(y_k,\hat{x}_k^{\mr{p},\Delta}) =0, \: k \in \hat{\mc{K}}.} \label{eq:nblcbenders0k}}
Doing so renders a new problem which is again amenable to nonlinear Benders decomposition and enjoys the same two advantages iterated previously.

\subsubsection{Reformulating the LBLC problem~\eqref{eqgroup:lblc} for decomposition}
Per contingency $k \in \hat{\mc{K}}$, introduce new decision variables $y_k$ with the same dimension as $\hat{x}_0^{\mr{p}}$. Remove constraint~\eqref{eq:lblc0k}.  Insert~\eqref{eq:nblcbenders0k} and the following constraint:
\eq{rCl}{\IEEEyesnumber 
	\hat{x}_0^{\mr{p}}&=&y_k, \:  k \in \hat{\mc{K}.} \label{eq:lblcbenderseq}}
	An advantage of the LBLC formulation is that since it is an LP,  there are convergence guarantees for Benders. A disadvantage, however, is that the incumbent optimal solution of the master problem at iteration $\ell$, denoted by $\hat{x}_0(\ell)$, is not immediately feasible for the base case. Therefore, per iteration, the recovery problem~\eqref{eqgroup:recovery} must be solved to calculate \texttt{solution1.txt} as a backup prior to the expiration of the time limit. 

\vspace{-2.5mm}
\section{Corrective SCOPF (``code 2") solution} \label{sec:code2heuristics}
\vspace{-1mm}
So far we have discussed approaches for code 1 which eventually output a base case feasible solution $x_0^*$ to \texttt{solution1.txt}.  The purpose of code 2 is to provide a solution $x_k^*$ for $k \in \mc{K}$. Notice here that $\mc{K}$ is the entire contingency set given by the data files and not the \emph{reduced} set discussed in Section~\ref{sec:code1heuristics}.  Ideally, given the obtained base case feasible solution $x_0^*$ by code 1,  code 2 is responsible for solving the following optimization problem:
\eq{rCl}{\IEEEyesnumber \label{eqgroup:correction} \IEEEyessubnumber*
		\minimize_{x_k,s_k} &&\frac{1}{K} \sum\limits_{k\in \mc{K}} c_k(s_k) \IEEEeqnarraynumspace \label{eq:correctionobj} \\
	\subjectto 
	f_k(x_k) &=& s_{k}^{f+}-s_{k}^{f-}, \: k \in \mc{K} \label{eq:correctionfk}\\
	g_k(x_k)&\le&s_{k}^{g}, \: k \in \mc{K}\label{eq:correctiongk}\\
	\IEEEeqnarraymulticol{3}{c}{h_k(x_0^{\mr{p},\mr{q},\mr{v}*},x_k^{\mr{p},\mr{q},\mr{v},\Delta}) =0, \: k \in \mc{K}} \label{eq:correction0k} \\
	\IEEEeqnarraymulticol{3}{c}{x_k \in \mc{X}_k, s_k \in \mc{S}_k, \: k \in \mc{K}. \label{eq:correctionboxes}}
}
Indeed, since $x_0^*$ is known in~\eqref{eq:correction0k}, code 2 offers a completely decentralized formulation per $k \in \mc{K}$. Nevertheless, as iterated previously, the structure of the nonconvexity of function $f_k(\cdot)$ and the nonconvexity and nondifferentiability of function $h_k(\cdot)$ remain as challenges for code 2. Due to the very stringent time limits imposed on code 2 by the GO competition (2 seconds average per contingency), we approximate~\eqref{eqgroup:correction} per contingency $k \in \mc{K}$:
\eq{rCl}{\IEEEyesnumber \label{eqgroup:correctionapprox} \IEEEyessubnumber*
	\minimize_{x_k,s_k} &&c_k(s_k) \IEEEeqnarraynumspace \label{eq:correctionapproxobj} \\
	\subjectto 
	f_k(x_k) &=& s_{k}^{f+}-s_{k}^{f-},  \label{eq:correctionapproxfk}\\
	g_k(x_k)&\le&s_{k}^{g}, \: k \in \mc{K}\label{eq:correctionapproxgk}\\
	\IEEEeqnarraymulticol{3}{c}{\hbar_k(x_0^{\mr{p},\mr{v}*},x_k^{\mr{p},\mr{v},\Delta}) =0, \: k \in \mc{K} \label{eq:correctionapprox0k}}\\
	\IEEEeqnarraymulticol{3}{c}{x_k \in \mc{X}_k, s_k \in \mc{S}_k. \label{eq:correctionapproxboxes}}
}
In~\eqref{eqgroup:correctionapprox}, the difficulty of handling~\eqref{eq:correction0k} is replaced by using its inner approximation~\eqref{eq:correctionapprox0k}. Notice that due to inner approximation, a feasible point of~\eqref{eq:correctionapprox0k}  is readily feasible for~\eqref{eq:correction0k}. The latter implies that if~\eqref{eqgroup:correctionapprox} is feasible, then the solution per $k\in \mc{K}$ is  feasible for the original code 2 problem~\eqref{eqgroup:correction}.  

\vspace{-3mm}
\section{Practical Implementation}
\label{sec:implementation}
\begin{figure}
    \centering
    \includegraphics[scale=0.35]{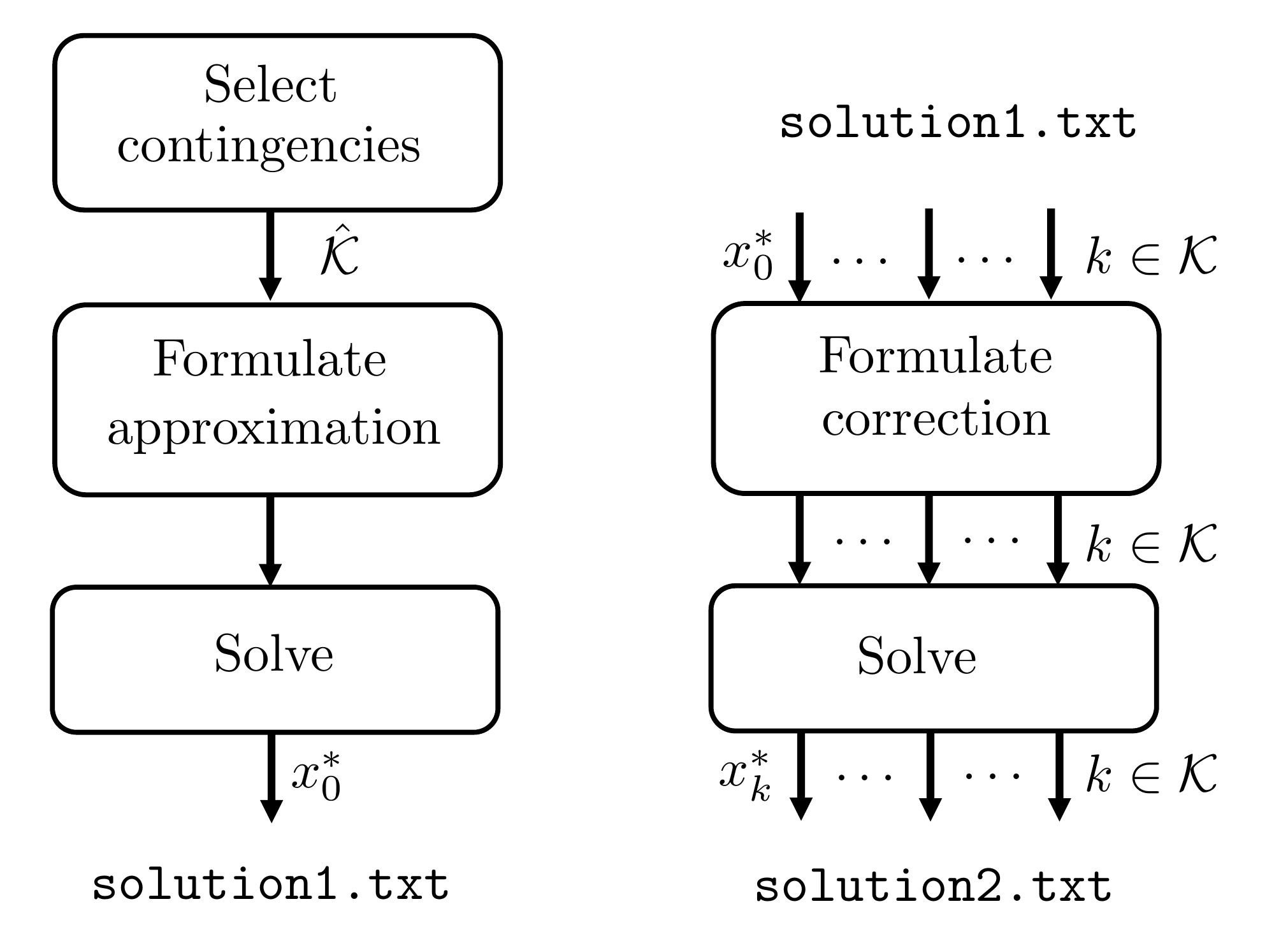}
    \caption{Sequential tasks of code 1 (left) and code 2 (right).  Multiple parallel arrows depict that the data transfer as well the task at the arrows' destinations may be performed in parallel.}
    \label{fig:overall}\vspace{-3mm}
\end{figure}

This section clarifies how the approaches in Sections~\ref{sec:code1heuristics} and~\ref{sec:code2heuristics} can be implemented on the multi-node competition evaluation platform~\cite{ARPAEPlatform}.  After reading data files for a specific test network, the crucial tasks of the submitted algorithms for code 1 and code 2 are shown on the left- and right-hand sides of Fig.~\ref{fig:overall}, respectively. 

To elaborate on Fig.~\ref{fig:overall}, code 1 first reduces the original contingency set $\mc{K}$ to the set $\hat{\mc{K}}$ following the methods in Section~\ref{sec:conselection}. Next, it chooses one of the approximate formulations, namely, NBNC~\eqref{eqgroup:nbnc}, NBLC~\eqref{eqgroup:nblc}, or LBLC~\eqref{eqgroup:lblc}. Finally, code 1 solves the selected approximate formulation either by passing its extensive form directly to a solver or by using Benders decomposition of Section~\ref{sec:Benders}.  Upon termination, code 1 returns \texttt{solution1.txt}.

Code 2 first reads \texttt{solution1.txt} and extracts the base case solution $x_0^*$. Then, it formulates the approximate correction problems~\eqref{eqgroup:correctionapprox} for all contingencies $k \in \mc{K}$. Upon solving these correction problems, the solutions $x_k^*$ for  all $k \in \mc{K}$ are printed in \texttt{solution2}.  The tasks of code 2 may be parallelized  across compute nodes available on the computer cluster---indicated by multiple parallel arrows in Fig.~\ref{fig:overall}. As the implementation of code 2 is relatively straightforward, in what follows, we focus solely on design differences for code 1 based on discussions of Section~\ref{sec:approxformulations} and~\ref{sec:Benders}.

\begin{figure}
\centering
\includegraphics[scale=0.35]{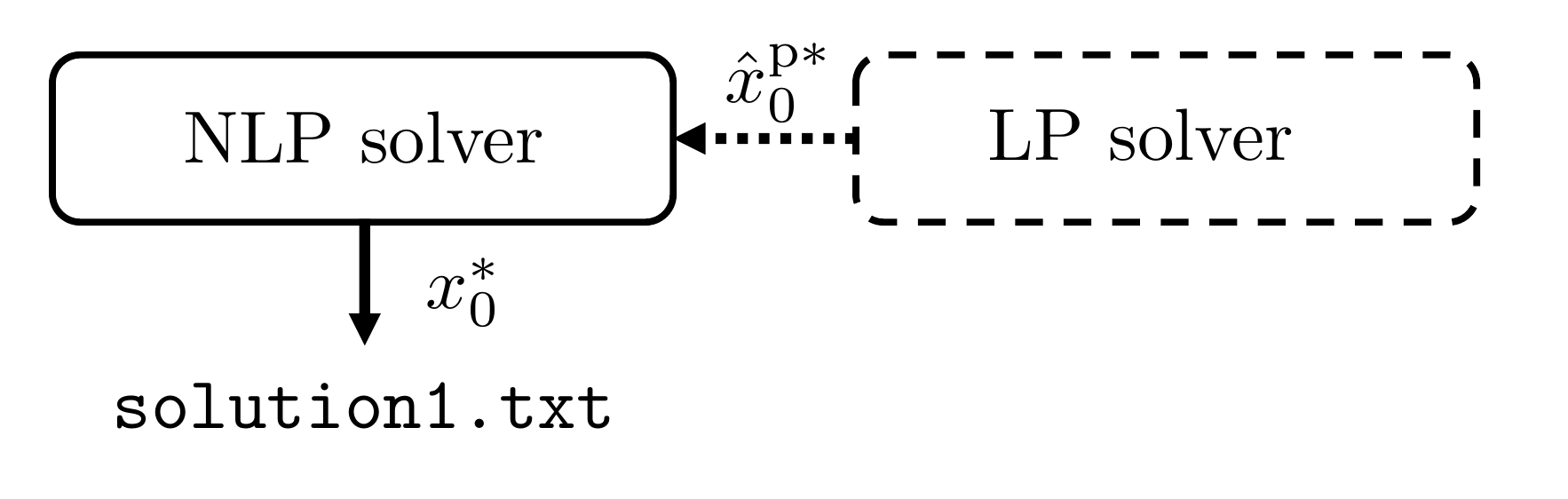}{}
    \caption{The solve task for code 1 when the extensive form of the approximate formulation is passed  directly to a solver.  For NBNC~\eqref{eqgroup:nbnc} and NBLC~\eqref{eqgroup:nblc} formulations, the dashed box is ignored and only an NLP solver is required to calculate \texttt{solution}. For the LBLC formulation~\eqref{eqgroup:lblc}, an LP solver is called first to solve for $x_0^{\mr{p}^*}$. An NLP solver then uses   $x_0^{\mr{p}^*}$ to solve the recovery problem~\eqref{eqgroup:recovery} to obtain $x_0^*$ and return \texttt{solution1.txt}. }\vspace{-3mm}
    \label{fig:lblc}
\end{figure}

\begin{figure}[t]
\centering
\includegraphics[scale=0.35]{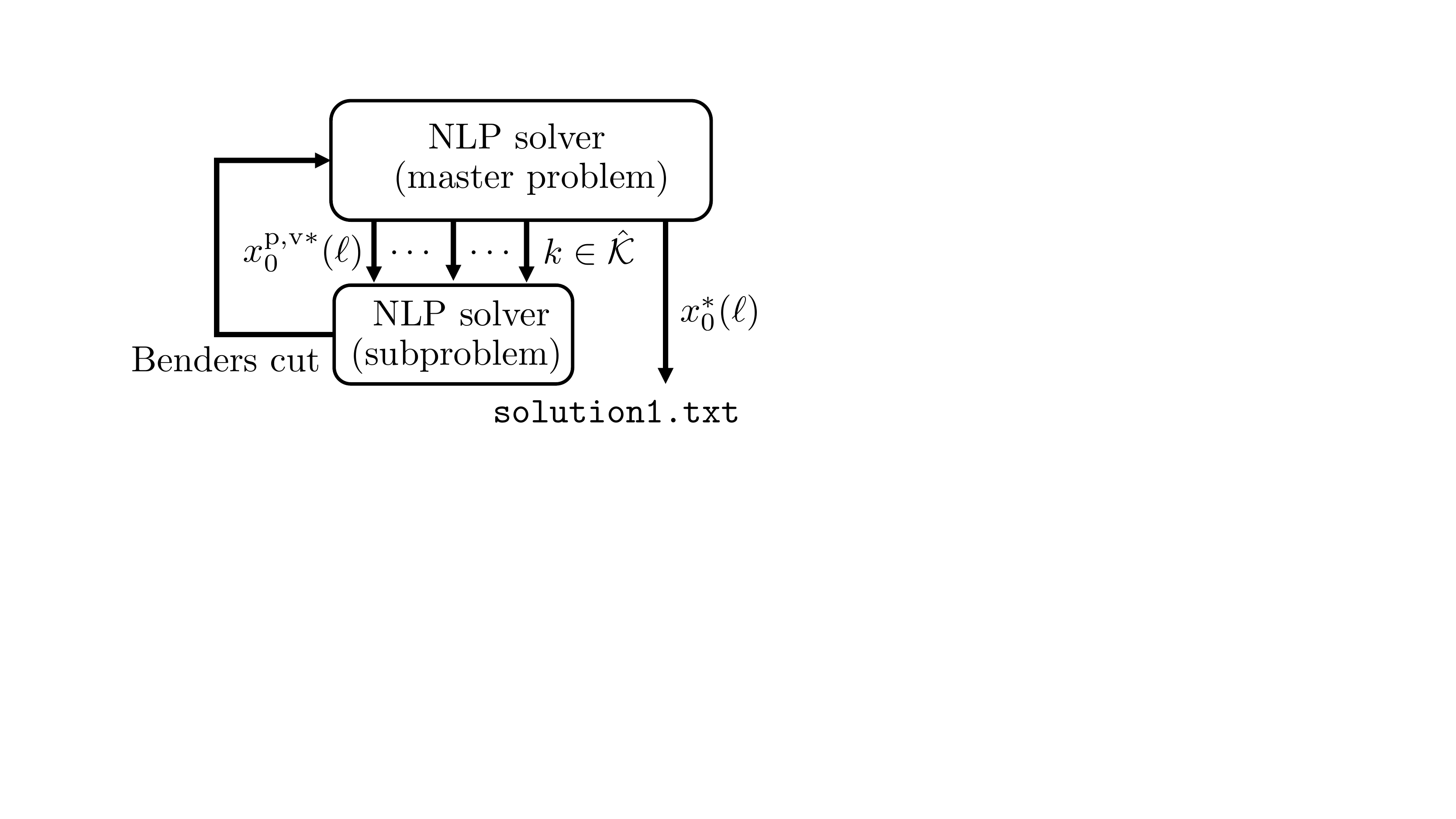}{}
    \caption{Iteration $\ell$ of the Benders decomposition  applied to NBNC.  Parallel arrows indicate that the data transfer and the tasks at the arrows' destination are conducted in parallel.  Here, \texttt{solution1.txt} is derived every iteration.}
    \label{fig:nbncb}\vspace{-5mm}
\end{figure}

 Figure~\ref{fig:lblc} shows subtasks of the code 1 solve task when the selected approximate formulation is passed  directly to a solver. Concretely, for the NBNC~\eqref{eqgroup:nbnc} and NBLC~\eqref{eqgroup:nblc} formulations, only an NLP solver is required to obtain the base case solution $x_0^*$ and the dashed box is ignored. For the LBLC formulation~\eqref{eqgroup:lblc},  an LP solver is called to obtain $\hat{x}_0^{\mr{p}*}$ prior to calling an NLP solver for the recovery problem~\eqref{eqgroup:recovery} to yield $x_0^*$. These methods output \texttt{solution1.txt} after the solver terminates.
In contrast to Fig.~\ref{fig:lblc}, one iteration of Benders applied to the NBNC formulation~\eqref{eqgroup:nbncbenders} is shown in Fig.~\ref{fig:nbncb}. In one compute node, at iteration $\ell$ the master problem is solved first using an NLP solver to obtain the  solution $x_0^*(\ell)$. This compute node writes  $x_0^*(\ell)$ to \texttt{solution1.txt} and communicates its entries $x_0^{\mr{p},\mr{v}*}(\ell)$ according to~\eqref{eq:nbncbenderseq}  to other compute nodes where subproblems pertaining to  specific contingencies $k \in \hat{\mc{K}}$ are solved in parallel.  Necessary information to form Benders cuts are passed back to the compute node containing the master problem. Notice a major difference between Fig.~\ref{fig:lblc} and Fig.~\ref{fig:nbncb} is that in the latter \texttt{solution1.txt} is written after every iteration. Thus, even if the algorithm does not converge, a feasible base case solution is still available.  Applying Benders to NBLC has the same configuration except that  subproblems are LPs and are solved faster than their nonlinear counterparts.

\begin{figure}[t]
\centering
\includegraphics[scale=0.35]{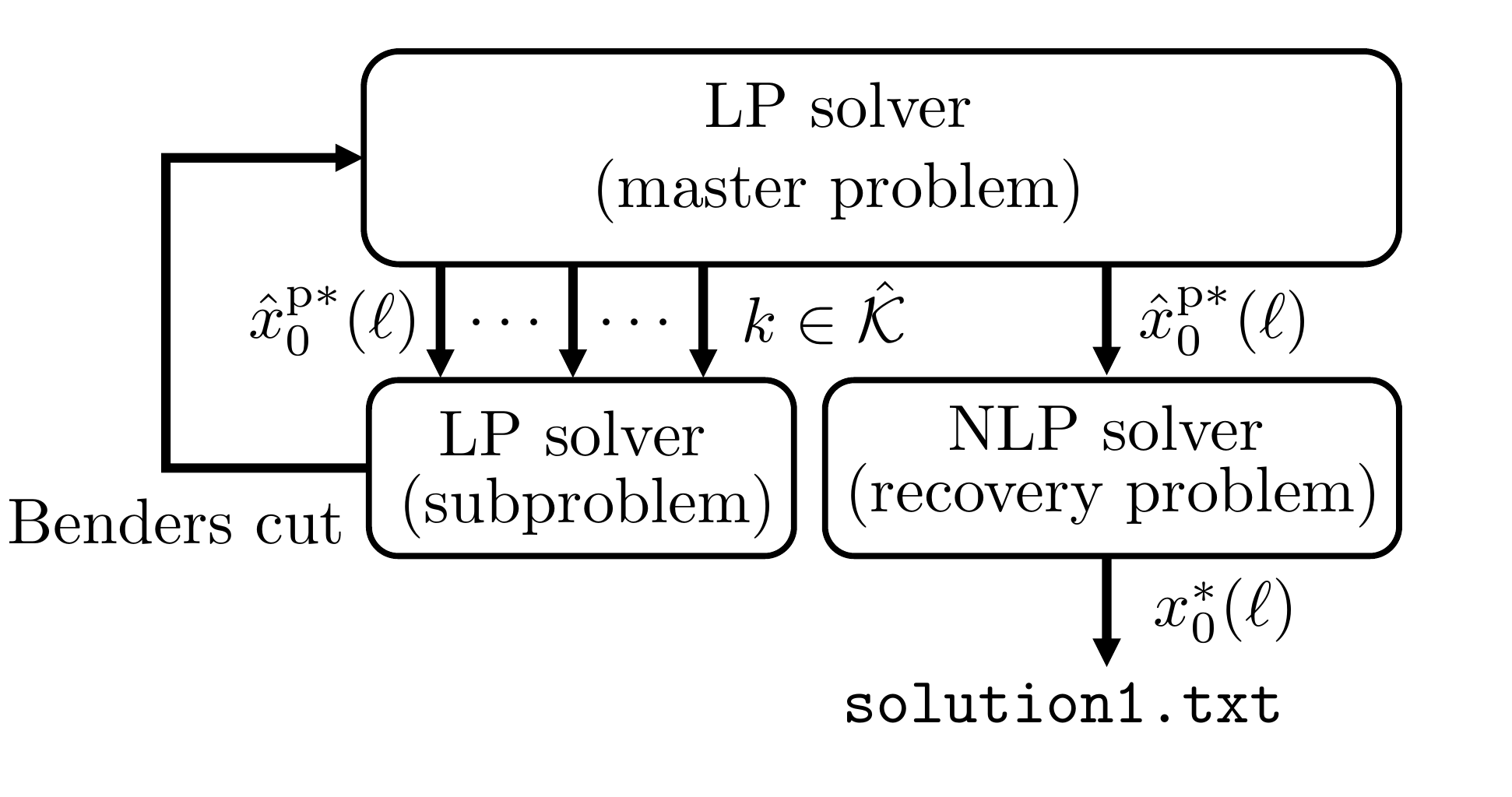}{}
    \caption{Design of iteration $\ell$ of the Benders decomposition  applied to LBLC. This method is different than applying Benders on NBNC and NBLC formulations  in that it requires solving the recovery problem~\eqref{eqgroup:recovery} in parallel with the contingency subproblems for $k \in \hat{\mc{K}}$.}    \label{fig:lblcb}\vspace{-3mm}
\end{figure}

One iteration of Benders decomposition method applied to  LBLC  is depicted in Fig.~\ref{fig:lblcb}. This design is  different than applying Benders to NBNC and NBLC  because the LBLC ~\eqref{eqgroup:lblc} is reliant on solving the recovery problem~\eqref{eqgroup:recovery} to ensure a feasible base case solution $x_0^*$.  Therefore, at every iteration $\ell$ of the Benders algorithm, the optimal solution of the master problem $\hat{x}_0^{\mr{p}*}(\ell)$ is not only passed to compute nodes that solve the LP subproblems, but it is also sent to a designated compute node that solves the recovery problem~\eqref{eqgroup:recovery}.   

\vspace{-2mm}
\section{Simulation Results}\label{sec:numericaltests}
The three NBNC, NBLC, and LBLC formulations, their respective Benders decompositions, together with the two contingency selection methods based on size and real-time implementation result in twelve distinct methodologies for solving the  SCOPF problem~\eqref{eqgroup:scopf}. In this section, we apply these twelve methods and compare their scores against one another. For demonstration purposes two test networks have been selected from 
ARPA-E GO data-sets~\cite{ARPAEDatasets}. The first  is  scenario 1 of \texttt{Network\textunderscore03R-10} which is a small network with 793 buses, 912 branches, and 91 contingencies. The second is scenario 1 of \texttt{Network\textunderscore07R-10}  which is a larger network with 2312 buses, 3013 branches, and 990 contingencies. 

\begin{figure}
\centering
\includegraphics[scale=0.3]{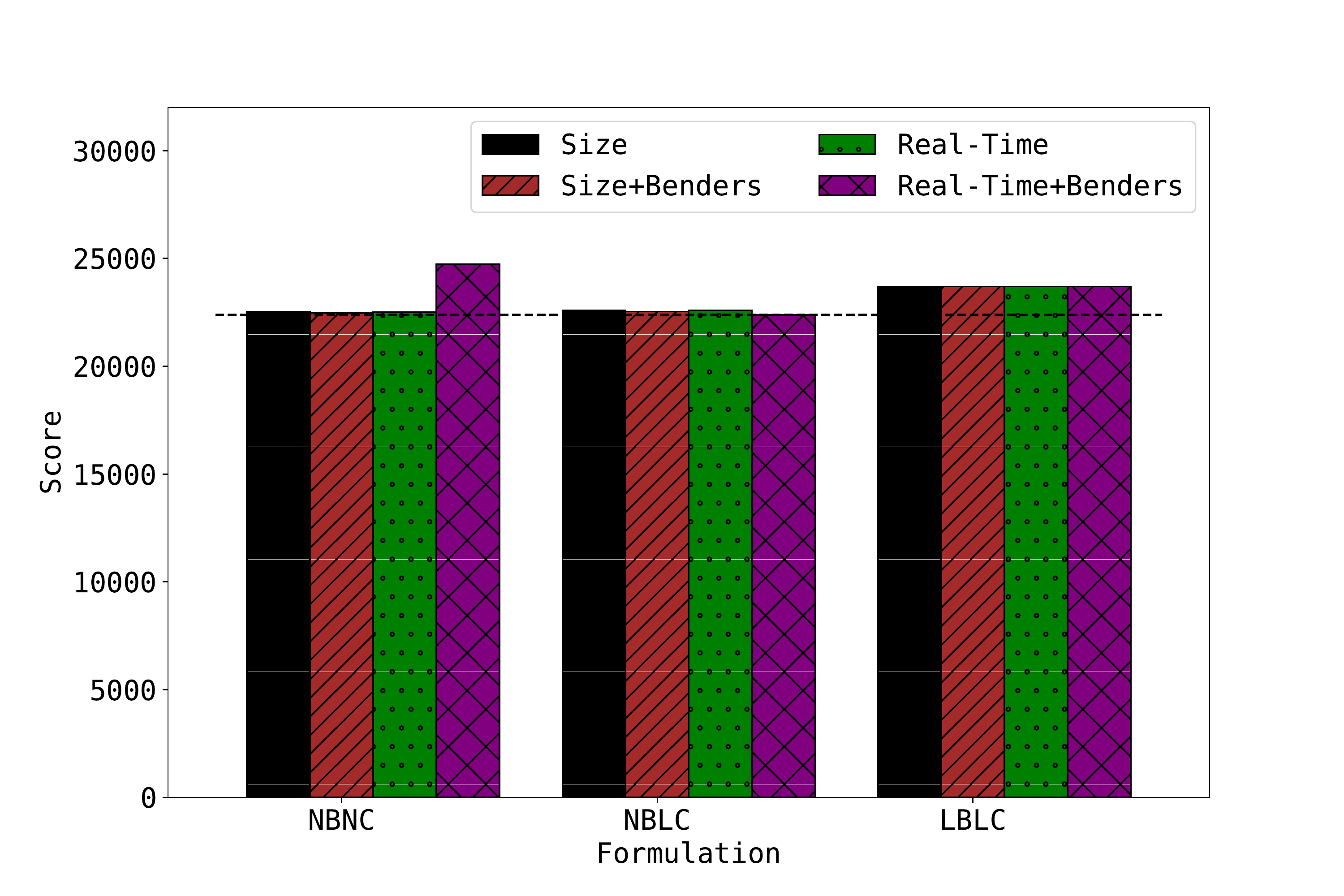}{}\vspace{-3mm}
    \caption{Results on a small network with 793 buses, 912 branches and 91 contingencies. All twelve heuristics achieve a relatively ``good" score that is much smaller than $c^{\mr{slack}}$.}    \label{fig:network3}\vspace{-3mm}
\end{figure}

For each of the twelve methods, we select $\hat{K}=5$. Code 1 and code 2 are set up based on the schematic in Fig.~\ref{fig:overall} and code 2 uses 144 processors. Upon termination of code 1 and code 2,  an evaluation algorithm analogous to Section~\ref{sec:evaluation} is run. 
The calculated scores for the smaller network are shown using a bar graph in Fig.~\ref{fig:network3}.   The x-axis represents the approximate formulation and the y-axis measures the scores.  The legends help distinguish between the size or real-time contingency selection methods and the application of Benders decomposition.  Figure~\ref{fig:network3} indicates that for small networks, all twelve methods achieve similar scores that are much lower than $c^{\mr{slack}}$.
The best solution among the twelve turns out to be using real-time contingency selection and Benders decomposition on the NBLC formulation, yielding a score of approximately 22380.433.

\begin{figure}
\centering
\includegraphics[scale=0.3]{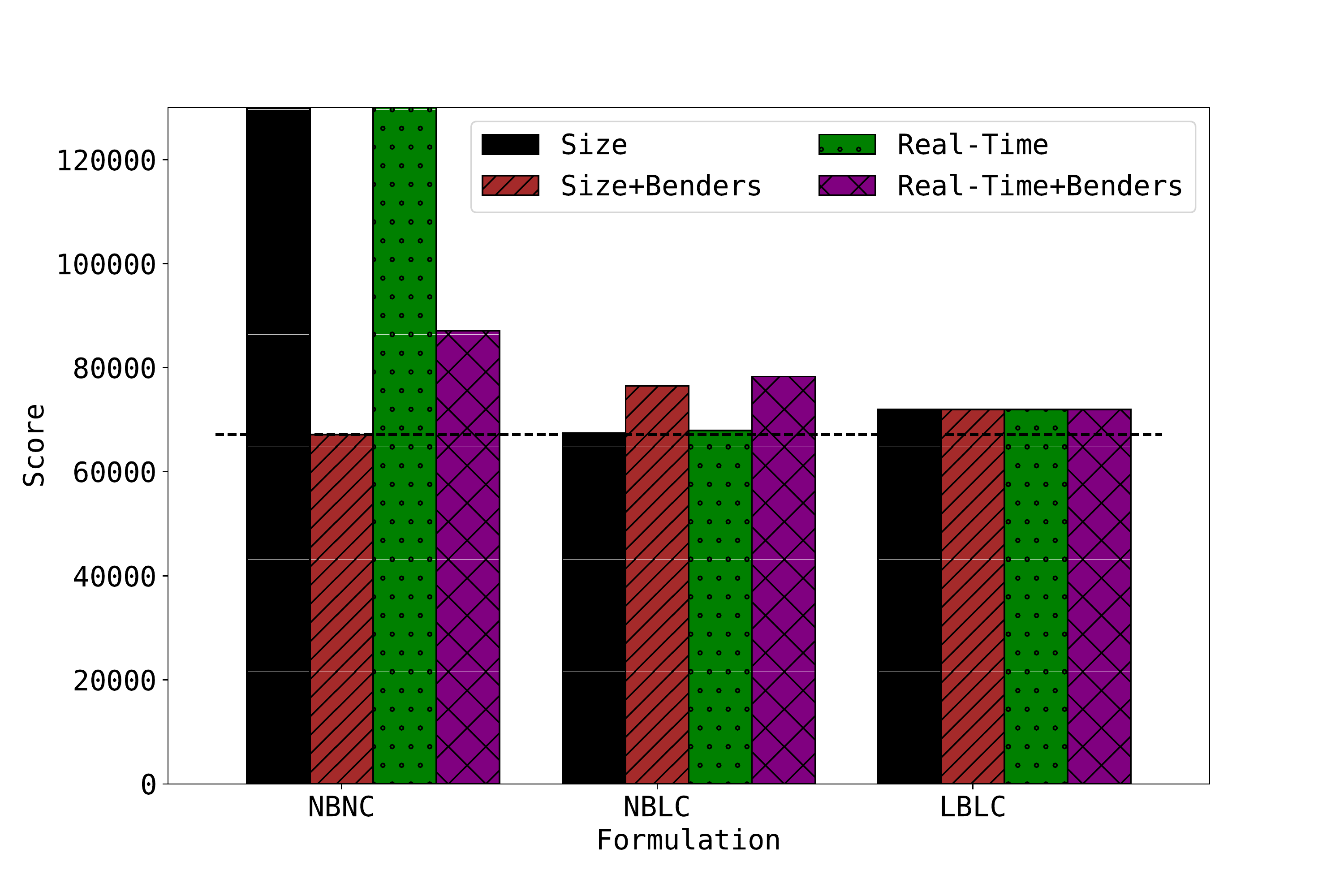}{}\vspace{-3mm}
    \caption{Results on a larger network with 2312 buses, 3013 branches and 990 contingencies. The extensive forms of the NBNC formulation time out and produce the worst-case solution.  The remaining heuristics obtain a ``good" score that is much smaller than $c^{\mr{slack}}$.}    \label{fig:network7}\vspace{-3mm}
\end{figure}

The calculated scores for the larger network are depicted in Fig.~\ref{fig:network7}. On this larger network, code 1 from extensive formulations of NBNC times out and thus outputs the worst-case base feasible solution. Code 2 run on the worst-case base feasible solution  does not finish in time also and violates the average 2 seconds per contingency time limit. Therefore, their corresponding scores are set to $c^{\mr{slack}}$ which exceed the y-axis in Fig.~\ref{fig:network7}. The remaining methods achieve much lower scores than $c^{\mr{slack}}$. For this network, the best score is 67158.175 and it comes from contingency selection based on size and applying Benders decompostion  to the NBNC formulation.

\vspace{-2mm}
\section{Conclusions}\label{sec:conclusion}
This paper presents computationally efficient solutions for solving comprehensive industry-based SCOPF formulation, which was thoroughly tested and performed well on a suite of both synthetic and actual industry networks. We proposed and compared twelve solution approaches to solve this mixed integer nonconvex stochastic problem and find the best feasible solution for large-scale systems in near \textit{real-time} (under ten minutes). Our proposed techniques are categorized in three groups contingency selection, formulation approximation, and decompositions. Future work includes extending our framework to include discrete transformer tap settings, switchable shunts, phase shifting transformers, price-responsive demand, generator ramp rates and transmission switching.
 \begin{figure}
 \centering
 \includegraphics[scale=0.3]{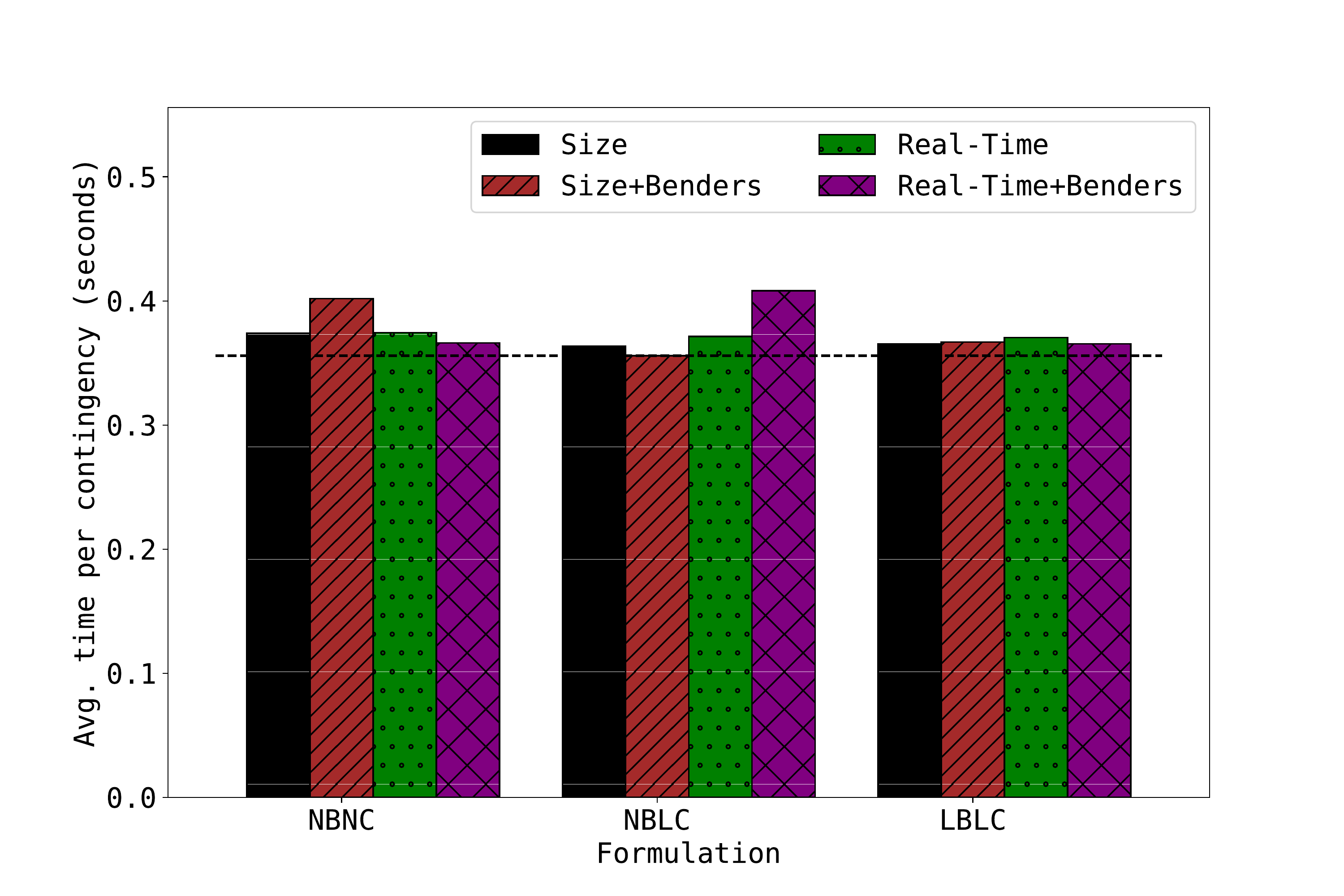}{}
     \caption{Average time per contingency on a small network.}    \label{fig:network3time}
 \end{figure}
 \begin{figure}
 \centering
 \includegraphics[scale=0.3]{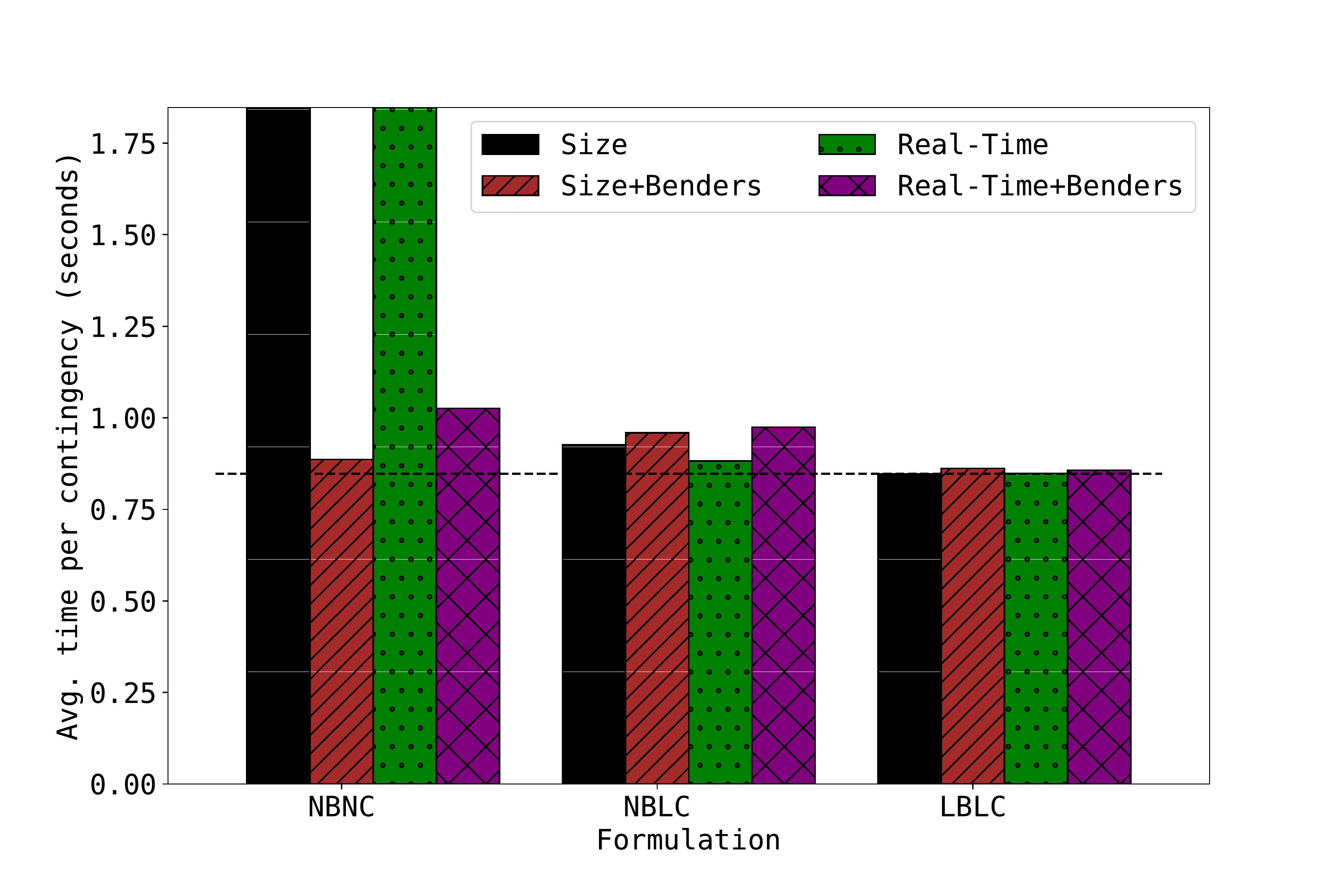}{}
     \caption{Average time per contingency on a larger network.}    \label{fig:network7time}
 \end{figure}

 
\bibliographystyle{IEEEtran}
\bibliography{biblio}

\end{document}